\RequirePackage{fix-cm}
\documentclass[twocolumn,10pt,a4paper]{article}          
\usepackage{graphicx}
\usepackage{amsmath}
\usepackage[authoryear,round]{natbib}
\usepackage{amsfonts,amssymb,amsthm}
\usepackage{bbm}
\usepackage{url}
\usepackage[cm]{fullpage}

\newtheorem{theorem}{Theorem}[section]

\newtheorem{proposition}[theorem]{Proposition}
\newtheorem{corollary}[theorem]{Corollary}

\newenvironment{remark}[1][Remark]{\begin{trivlist}
\item[\hskip \labelsep {\bfseries #1}]}{\end{trivlist}}

\newcommand{\Rset}{\mathbb{R}}								
\newcommand{\Xset}{\mathbb{X}}

\newcommand{\prob}{{\mathbb{P}}}
\newcommand{\x}{\mathbf{x}}
\newcommand{\X}{\mathbf{X}}
\newcommand{\y}{\mathbf{y}}
\newcommand{\Y}{\mathbf{Y}}
\newcommand{\+}{_{n+1}}

\newcommand{\An}{\mathcal{A}_n}
\newcommand{\nlaw}{\mathcal{N}}						
\newcommand{\esp}{\mathbb{E}}	

\newcommand{ \ykip} {y^{(k)}_{i+}}
\newcommand{ \ykim} {y^{(k)}_{i-}}
\newcommand{ \ykjp} {y^{(k)}_{j+}}
\newcommand{ \ykjm }{y^{(k)}_{j-}}
\newcommand{\one}{\mathbbmss{1}}

\begin{document}

\title{Multiobjective optimization using Gaussian process emulators via stepwise uncertainty reduction}
\author{Victor Picheny\\
INRA, 31326 Castanet Tolosan, France\\
Tel.: +33-5-61 28 54 39\\
victor.picheny@toulouse.inra.fr
}


\maketitle

\begin{abstract}
Optimization of expensive computer models with the help of Gaussian process emulators in now commonplace. 
However, when several (competing) objectives are considered, choosing an appropriate sampling strategy remains an open question.
We present here a new algorithm based on stepwise uncertainty reduction principles to address this issue. 
Optimization is seen as a sequential reduction of the volume of the excursion sets below the current best solutions, 
and our sampling strategy chooses the points that give the highest expected reduction. 
Closed-form formulae are provided to compute the sampling criterion, avoiding the use of cumbersome simulations.
We test our method on numerical examples, showing that it provides an efficient trade-off between exploration and intensification.

\textbf{keywords} Kriging; EGO; Pareto front; Excursion sets
\end{abstract}
\section{Introduction}
We consider the problem of simultaneous optimization of several objective functions over a design region $\Xset \subset \mathbb{R}^d$:
\begin{eqnarray}
 \min & y^{(1)}(\x), \ldots, y^{(q)}(\x),   \nonumber
\end{eqnarray}
where $y^{(i)}: \Xset \rightarrow \mathbb{R}$ are outputs of a complex computer code. The objectives being typically conflicting, there exists no unique minimizer, and the goal is to identify the set of optimal solutions, called Pareto front \citep{collette2003multiobjective}. Defining that a point dominates another if \textit{all} his objectives are better, the Pareto front $\mathbb{X}^*$ is the subset of the non-dominated points in $\Xset$:
\begin{eqnarray}
 \forall \x^* \in \mathbb{X}^*, \forall \x \in \mathbb{X}, \exists k \in \{1, \ldots, q\} \text{ such that } \nonumber\\
  y^{(k)}(\x^*) \leq y^{(k)}(\x). \nonumber
\end{eqnarray}
When the computational cost of a single model evaluation is high, a well-established practice consists of using Gaussian process (GP) emulators to approximate the model outputs and guide the optimization process. Following the seminal article of \citet{jones1998efficient} and its Efficient Global Optimization (EGO) algorithm for single objective optimization, several strategies have been proposed in the past few years to address the multi-objective problem \citep{knowles2006parego,keane2006statistical,ponweiser2008multiobjective,wagner2010expected}. 
They consist in evaluating sequentially the computer model at the set of inputs that maximizes a so-called \textit{infill criterion}, derived from the GP emulator, that expresses a trade-off between exploration of unsampled areas and sampling intensification in promising regions. While the single objective case has been extensively discussed \citep{jones2001taxonomy,wang2007review}, finding efficient and statistically consistent infill criteria for the multi-objective case remains an open question.

Alternatively to the EGO paradigm, \textit{stepwise uncertainty reduction} (SUR) strategies aim at reducing, by sequential sampling, an uncertainty measure about a quantity of interest. In a single objective optimization context, \citet{villemonteix2009informational} defined the Shannon entropy of the maximizer (computed using the GP model) as an uncertainty measure: a smaller entropy implies that the maximizer is well-identified. They show that their approach outperforms the EGO strategy on a series of problems. Another example in a reliability assessment context can be found in \citet{bect2012sequential}. In general, SUR approaches allow to define policies rigorously with respect to a given objective, resulting in very good performances. However, they are often challenging to use in practice, as they rely on very expensive GP simulations. 

We propose here a novel SUR strategy to address the multi-objective problem. It is based on a measure of uncertainty of the current identification of the Pareto front $\mathbb{X}^*$, hence avoiding some of the drawbacks of the existing criteria (hierarchy between objectives, difficult-to-tune parameters, etc.). Following \citet{chevalier2012fast}, explicit formulae for the expected uncertainty reduction are provided, avoiding the need to rely on simulations.

The paper is organized as follows: section \ref{sec:background} presents the GP model and the basics of GP-based optimization. Then, we describe our SUR strategy for a single objective in section \ref{sec:unconstrained} and for several objectives in section \ref{sec:multi}. We provide some numerical experiments in section \ref{sec:illustration} and compare our method to the state-of-the-art. Finally, advantages and drawbacks of the method are discussed in section \ref{sec:discussion}.
\section{Some concepts of Gaussian-process-based optimization}\label{sec:background}
\subsection{Gaussian process emulation}
We consider first the emulation of a single computer response $y$. The response is modelled as
\begin{equation}
 Y(.) = \mathbf{f}(.)^{T}\boldsymbol{\beta} + Z(.),
\end{equation}
where $\mathbf{f}(.)^{T} = \left( f_1(.), \ldots, f_p(.) \right)$ is a vector of trend functions, $\boldsymbol{\beta}$ a vector of (unknown) coefficients and $Z(.)$
is a Gaussian process (GP) $Z$ with zero mean and known covariance kernel $k$ \citep{cressie1993statistics,rasmussen2006gaussian}. 
Let us call $\An$ the event: $$\left\{ Y(\mathbf{x}_1)=y_1, \ldots, Y(\mathbf{x}_n)=y_n \right\};$$
conditionally on $\An$, the mean and covariance of $Y$ are given by:
\begin{eqnarray}
m_{n}(\mathbf{x})&=& \esp \left( Y(\mathbf{x}) \big| \An  \right) =
 \nonumber \\
 &=&
\mathbf{f}(\mathbf{x})^T \boldsymbol{\hat  \beta}
+\mathbf{k}_{n}(\mathbf{x})^{T}\mathbf{ K}_{n}^{-1}({\mathbf{y}}_{n}- \mathbf{F}_n\boldsymbol{\hat  \beta}), \nonumber \\ 
c_{n}(\mathbf{x},\mathbf{x}') &=& cov \left( Y(\mathbf{x}), Y(\mathbf{x}') \big| \An  \right)  \nonumber \\
&=&
k(\mathbf{x}, \mathbf{x}')-\mathbf{k}_{n}(\mathbf{x})^{T}\mathbf{ K}_{n}^{-1}\mathbf{k}_{n}(\mathbf{x}') \nonumber \\
&+&
\left(\mathbf{f}(\mathbf{x})^T - \mathbf{k}_{n}(\mathbf{x})^{T}\mathbf{ K}_{n}^{-1} \mathbf{F}_n \right)^T
\left( \mathbf{F}_{n}^{T}\mathbf{ K}_{n}^{-1}\mathbf{F}_{n} \right)^{-1} \nonumber \\
& &
\left(\mathbf{f}(\mathbf{x}')^T - \mathbf{k}_{n}(\mathbf{x}')^{T}\mathbf{ K}_{n}^{-1} \mathbf{F}_n \right), \nonumber 
\end{eqnarray}
where 
\begin{itemize}
 \item ${\mathbf{y}}_{n}=\left( {y}_{1},\ldots, {y}_{n} \right)^{T}$ are the observations,
 \item $\mathbf{K}_{n}=\left( k(\mathbf{x}_{i}, \mathbf{x}_{j}) \right)_{1\leq i,j \leq n}$ is the observation covariance matrix,
 \item  $\mathbf{k}_{n}(\mathbf{x})^{T}=\left( k(\mathbf{x}, \mathbf{x}_{1}), \ldots, k(\mathbf{x}, \mathbf{x}_{n}) \right)$,
 \item  $\mathbf{F}_n = \left( \mathbf{f}(\mathbf{x}_1)^T, \ldots,  \mathbf{f}(\mathbf{x}_n)^T \right)^T$, and
 \item $\boldsymbol{\hat  \beta} = \left( \mathbf{F}_{n}^{T}\mathbf{ K}_{n}^{-1}\mathbf{F}_{n} \right)^{-1} \mathbf{F}_{n}^{T}\mathbf{
K}_{n}^{-1} {\mathbf{y}}_{n}$ is the best linear unbiased estimate of $\boldsymbol{\beta}$. 
\end{itemize}
In addition, the \textit{prediction variance} is defined as
  $$s_n^2(\mathbf{x}) = c_{n}(\mathbf{x},\mathbf{x}).$$ 

The covariance kernel depends on parameters that are usually unknown and must be estimated from an initial set of responses. Typically, maximum likelihood estimates are obtained by numerical optimization and used as face value, the estimates being updated when new observations are added to the model. The reader can refer to \citet{stein1999interpolation} (chapter 6), \citet{rasmussen2006gaussian} (chapter 5) or \citet{roustant2012dicekriging} for detailed calculations and implementation issues.

When several functions $y^{(1)}, \ldots, y^{(q)}$ are predicted simultaneously, it is possible to take their dependency into account \citep{kennedy2001bayesian,craig2001bayesian}. However, in this work we consider all the processes $Y^{(i)}$ independent, hence modelled as above, which is in line with current practice.

\subsection{Gaussian-process-based optimization with a single objective}
The EGO strategy, as well as most of its modifications, is based on the following scheme. An initial set of observations is generated, from which the GP model is constructed and validated. 
Then, new observations are obtained sequentially, at the point in the design space that maximizes the infill criterion, and the model is updated every time a new observation is added to the training set. The two later steps are repeated until a stopping criterion is met.

The expected improvement criterion ($EI$) used in EGO relies on the idea that progress is achieved by performing an evaluation at step $n$ if the $(n+1)^{\text{th}}$ design has a lower objective function value than any of the $n$ previous designs. Hence, the \textit{improvement} is defined as the difference between the current observed minimum and the new function value if it is positive, or zero otherwise, and $EI$ is its conditional expectation under the GP model:
\begin{equation}
EI(\mathbf{x}) = \mathbb{E}\left[ \max \left( 0, y_n^{\min} - Y\left(\mathbf{x} \right) \right) | \An \right],\nonumber
\end{equation}
where $y_n^{\min}$ denotes the current minimum of $y$ found at step $n$: $y_n^{\min} = \min (y_1, \ldots, y_n)$.

EGO is the one-step optimal strategy (in expectation) regarding improvement: at step $n$, the new measurement is chosen as 
\begin{equation}
 \x\+ = \arg \max_{\x \in \Xset} EI(\x), \nonumber
\end{equation}
which is in practice done by running an optimization algorithm. It has been shown in \citet{jones2001taxonomy} that EGO provides, among numerous alternatives, an efficient 
solution for global optimization.

\subsection{Gaussian-process-based optimization with several objectives}
Several adaptations of EGO to the multi-objective framework have been proposed; a review can be found in \citet{ponweiser2008multiobjective}. The main difficulty is that the concept of \textit{improvement} cannot be transferred directly, as the current best point is here a set, and the gain is measured on several objectives simultaneously. In \citet{knowles2006parego}, the objectives are aggregated in a single function using random weights, which allows using the standard EGO. \citet{keane2006statistical} derived an EI with respect to multiple objectives. \citet{ponweiser2008multiobjective} proposed an hypervolume-based infill criterion, where the improvement is measured in terms of hypervolume increase.

\section{Single objective optimization by stepwise uncertainty reduction}\label{sec:unconstrained}
We consider first the case of a problem with a single objective $y$ to minimize. In this section, 
we propose a new strategy in a form similar to EGO that uses an alternative infill criterion based on stepwise uncertainty reduction principles.
The adaptation of this criterion to the multi-objective case is presented in Section \ref{sec:multi}.

\subsection{Definition of an uncertainty measure for optimization}
The EGO strategy focuses on progress in terms of objective function value. It does not account (or only indirectly) for the knowledge improvement that a new measurement would provide to the GP model, 
nor for the discrepancy between the location of the current best design found and the actual minimizer (which is actually most users' objective).

Alternative sampling criteria have been proposed to account for these two aspects. In \citet{villemonteix2009informational}, the IAGO stategy chooses the point that minimizes the posterior Shannon entropy of the minimizer: the interest of performing a new observation is measured in gain of information about the location of the minimizer. Unfortunately, it relies on expensive GP simulations, which makes its use challenging in practice. \citet{gramacy2011optimization} proposed an \textit{integrated expected conditional improvement} to measure a global informational gain of an observation. In the noisy case, \citet{scott2011correlated} proposed a somehow similar \textit{knowledge gradient policy} that also measures global information gain. However, as both criteria rely on notions of improvement, it makes them difficult to adapt to the multiobjective case. The criterion we propose below address this issue.

Consider that $n$ measurements have been performed. As a measure of performance regarding the optimization problem, we consider the expected volume of excursion set below the current minimum $y_n^{\min}$:
\begin{equation}
 ev_n = \esp_{\Xset} \left[ \mathbb{P} \left( Y(\x) \leq y_n^{\min} | \An \right)  \right]. \label{eq:def1}
\end{equation}
Similarly to the Shannon entropy measure in IAGO, a large volume indicates that the optimum is not yet precisely located (see Figure \ref{fig:exSUR1}); on the contrary, a small volume indicates that very little can be gained by pursuing the optimization process. 
Following the stepwise uncertainty reduction paradigm, this volume is an uncertainty measure related to our objective (finding the minimizer of $y$); minimizing the uncertainty amounts to solving the optimization problem.

The probability $p_n(\x, y_n^{\min}) := \mathbb{P} \left( Y(\x) \leq y_n^{\min} | \An \right)$, which is often referred to as \textit{probability of improvement} \citep{jones2001taxonomy}, can be expressed in closed form, and Eq.~(\ref{eq:def1}) writes:
\begin{equation}
 ev_n =  \int_\Xset p_n(\x,y_n^{\min}) d\x = \int_\Xset \Phi \left( \frac{y_n^{\min} - m_n(\x)}{s_n(\x)} \right) d\x,
\end{equation}
where $\Phi(.)$ is the cumulative distribution function (CDF) of the standard Gaussian distribution.  
Hypothesizing that a measurement $y\+$ is performed at a point $\x\+$, its benefit can be measured by the reduction of the expected volume of excursion set
$\Delta = ev_n - ev\+ $, with:
\begin{eqnarray}
 ev\+ &=& \int_\Xset p\+(\x,\min\left( y_n^{\min}, y\+\right)) d\x \nonumber \\
      &=& \int_\Xset \Phi \left( \frac{  \min\left( y_n^{\min}, y\+\right) -  m\+(\x)}{s\+(\x)} \right) d\x.  \nonumber
\end{eqnarray}
Of course, $ev\+$ cannot be known exactly without evaluating $y\+$. However, we show in the following that its expectation can be calculated in closed form, leading to a suitable infill criterion. 
To do so, we first formulate a series of propositions in the next subsection.

\subsection{Probabilities updates}\label{sec:update}
An interesting property of the GP model is that, when a new observation $y\+ = y(\x\+)$ is added to the training set, its new predictive distribution can be expressed simply as a function of the old one \citep{emery2009kriging}:
\begin{eqnarray}
 m_{n+1}(\x) &=& m_n(\x) + \frac{c_n(\x,\x_{n+1})}{c_n(\x_{n+1},\x_{n+1})}\left({y}_{n+1} - m_n\left( \x_{n+1} \right)  \right);\nonumber\\ 
 s_{n+1}^2(\x) &=& s_n^2(\x) - \frac{c_n(\x,\x_{n+1})^2}{s_n^2(\x_{n+1})}. \label{eq:updatedsn2}
\end{eqnarray}
Note that only $m_{n+1}(\x)$ depends on the value of the new observation ${y}_{n+1}$. Now, conditionally on the $n$ first observations, $Y\+$ is a random variable (as the new observation has not yet been performed) with its moments given by the GP model: $$Y\+ \sim \nlaw \left( m_n(\x\+), s_n^2(\x\+) \right).$$ We can then define the future expectation $M_{n+1}(\x)$ (or any quantity depending on it) as a random variable conditionally on $\An$ and on the fact that the next observation will be at $\x_{n+1}$. 
This applies to any quantity depending on $Y\+$ or $M_{n+1}(\x)$, for instance, the probability of being below a threshold $a \in \Rset$:
\begin{equation}
 P\+(\x,a) = 
 \Phi \left( \frac{a - M\+(\x)}{s\+(\x)} \right). \nonumber
\end{equation}

\begin{proposition} \label{prop:4}
Without any restriction on the value of $Y\+$, the expectation of the future probability of being below the threshold is equal to the current probability:
\begin{eqnarray}
  \mathbb{P} \left(Y(\x) \leq a | \An, Y(\x\+)=Y\+ \right) &=& \mathbb{E} \left[ P\+(\x,a) \big| \An \right] \nonumber   \\
  &=& p_n (\x,a).\nonumber
\end{eqnarray}
\end{proposition}

\begin{proposition}\label{prop:3}
Conditioning further by $Y\+ \leq b$, the probability expectation writes in simple form using the Gaussian bivariate CDF:
\begin{eqnarray}
 q(\x, b, a) &:=& \mathbb{P} \left[ Y(\x) \leq a \big| \An, Y(\x\+)=Y\+, Y\+ \leq b \right] \nonumber \\
 &\times&\mathbb{P} \left[ Y\+ \leq b \big| \An \right] \nonumber \\
             &=&   \mathbb{E} \left[ P\+(\x,a) \times \one_{ Y\+ \leq b } \big| \An \right] \nonumber \\
             &=& \boldsymbol{\Phi}_{\rho} \left(\bar b, \tilde a \right),  \label{eq:prop3}
\end{eqnarray}
where $\boldsymbol{\Phi}_{\rho}$ is the Gaussian bivariate CDF with zero mean and covariance 
$\left[ \begin{matrix}
1  &  \rho\\
\rho  & 1   \end{matrix} \right]$,
$\bar b = \frac{b - m_n(\x\+) }{s_n(\x\+)}$, 
$\tilde a = \frac{a - m_n(\x)}{s_n(\x)}$ and 
$\rho = \frac{c_n(\x,\x\+)}{s_n(\x\+)s_n(\x)}.$
\end{proposition}

\begin{corollary}\label{corr:1} 
Similarly, conditioning by $Y\+ \geq b$ leads to:
\begin{eqnarray}
 r(\x, b, a) &:=& \mathbb{P} \big[ Y(\x) \leq a \big| \An, Y(\x\+)=Y\+, Y\+ \geq b \big] \nonumber \\
 &\times&\mathbb{P} \big[ Y\+ \geq b \big| \An \big]\nonumber \\
                     &=& \boldsymbol{\Phi}_{-\rho} \left(-\bar b, \tilde a \right). \label{eq:corr1}
\end{eqnarray}
\end{corollary}

The final proposition resembles Proposition \ref{prop:3}, but the fixed threshold $a$ is here replaced by $Y\+$:
\begin{proposition} \label{prop:5}
The expectation of the probability that $Y(\x)$ is smaller than $Y\+$, conditionally on $Y\+ \leq b$, is given by:
\begin{eqnarray}
  h(\x, b) &:=& \mathbb{P} \big[Y(\x) \leq Y\+ | \An, Y(\x\+)\nonumber \\
  &=&Y\+, Y\+ \leq b \big] \mathbb{P} \big[  Y\+ \leq b \big| \An \big] \nonumber \\
           &=&  \mathbb{E} \left[ P\+(\x,Y\+) \times \one_{ Y\+ \leq b } \big| \An \right] \nonumber \\
           &=&  \boldsymbol{\Phi}_{\nu} \left(\bar b, \eta \right), \label{eq:prop5}
\end{eqnarray}
with:
\begin{eqnarray}
\eta &=& \frac{m_n(\x\+) - m_n(\x)}{\sqrt{s_n^2(\x) + s_n^2(\x\+) - 2 c_n(\x,\x\+)}} \text{ and } \nonumber \\
\nu &=& \frac{c_n(\x,\x\+) - s_n^2(\x\+)}{s_n(\x\+)\sqrt{s_n^2(\x) + s_n^2(\x\+) - 2 c_n(\x,\x\+)}}. \nonumber
\end{eqnarray}
\end{proposition}

All the proofs are reported in Appendix \ref{sec:appendix2}.

\subsection{A Stepwise uncertainty reduction criterion}
Coming back to the SUR criterion, at step $n$ the future volume of excursion set $EV\+$ is a random variable, and its expectation is:
\begin{eqnarray} \label{eq:def3}
 EEV(\x\+) &:=&  \esp \left(EV\+ \big| \An, Y(\x\+) = Y\+\right) \nonumber \\
           &=&  \int_\Xset \esp \Big[ \Phi \left( \frac{  \min\left( y_n^{\min}, Y\+\right) -  M\+(\x)}{s\+(\x)} \right) \nonumber \\
           &&\big| \An, Y(\x\+) = Y\+ \Big] d\x.\nonumber
\end{eqnarray} 
Let $\varphi\left( y\+ \right)$ be the probability density function (PDF) of $Y\+$ conditionally on $\An$. We have:
\begin{eqnarray}
 && EEV(\x\+) \nonumber \\
 &   =& \int_\Xset \int_\Rset \Phi \left( \frac{  \min\left( y_n^{\min}, y\+\right) -  m\+(\x)}{s\+(\x)} \right) d\varphi\left( y\+ \right)d\x  \nonumber \\ 
 &   =& \int_\Xset \big[ \int_{-\infty}^{y_n^{\min}} \Phi \left( \frac{ y\+ -  m\+(\x)}{s\+(\x)} \right) \nonumber \\ 
 &   +&  \int_{y_n^{\min}}^{+\infty} \Phi \left( \frac{ y_n^{\min} -  m\+(\x)}{s\+(\x)} \right) d\varphi\left( y\+ \right) \big] d\x  \nonumber \\
 &   =& \int_\Xset \left[   h(\x,y_n^{\min}) + r(\x, y_n^{\min}, y_n^{\min}) \right] d\x.  \nonumber
\end{eqnarray}
The first term of the integrand is given by Eq.~(\ref{eq:prop5}) in Proposition \ref{prop:5}, with $b=y_n^{\min}$, and the second term is given by Eq.~(\ref{eq:corr1}) in Corrolary \ref{corr:1}, with $a=b=y_n^{\min}$, hence:
\begin{equation} \label{eq:unconscrit}
EEV(\x\+) = \int_\Xset{ \left[ \boldsymbol{\Phi}_{\nu}\left(\overline y_n^{\min}, \eta \right) + \boldsymbol{\Phi}_{\rho}\left( -\overline y_n^{\min}, \widetilde y_n^{\min} \right) \right] d\x },
\end{equation}
with:
$$\overline y_n^{\min} = {(y_n^{\min} - m_n(\x\+) )}/{s_n(\x\+)}$$
and $$\widetilde y_n^{\min} = {(y_n^{\min} - m_n(\x))}/{s_n(\x)}.$$

The SUR optimization strategy consists in adding the experiment that minimizes the expected volume of excursion set (or maximizes the difference), 
that is, the one-step optimal policy in terms of reduction of the uncertainty on the objective function minimizer:
\begin{equation}
 \x\+ = \arg \min_{\x^+ \in \Xset} EEV(\x^+) \label{eq:SURstrategy}
\end{equation}

\begin{remark}
In general, the probability of improvement $p_n(\x, y_n^{\min})$ is high where the prediction mean $m_n(.)$ is low and/or the prediction variance $s_n^2(.)$ is high. 
Simply choosing points that maximize $p_n(\x, y_n^{\min})$ is known to be inefficient \citep{jones2001taxonomy}, 
as it does not consider the amplitude of the gain in the objective function.
Here, $EEV(\x^+)$ strongly depends on the potential gain amplitude. 
Indeed, minimizing the expected volume relies on two mecanisms: reducing the local uncertainty and lowering the current minimum value ($y_n^{\min}$). The first is achieved by adding measurements in unsampled regions (high $s_n^2(.)$), the second in regions where this potential reduction is high.
Hence, the $EEV$ criterion can be seen as a mixed measure of uncertainty on the current minimum location and of potential gain in the objective function. 
\end{remark}

\subsection{Illustration}
\begin{figure*}[ht]
	\centering
	\includegraphics[height=55mm]{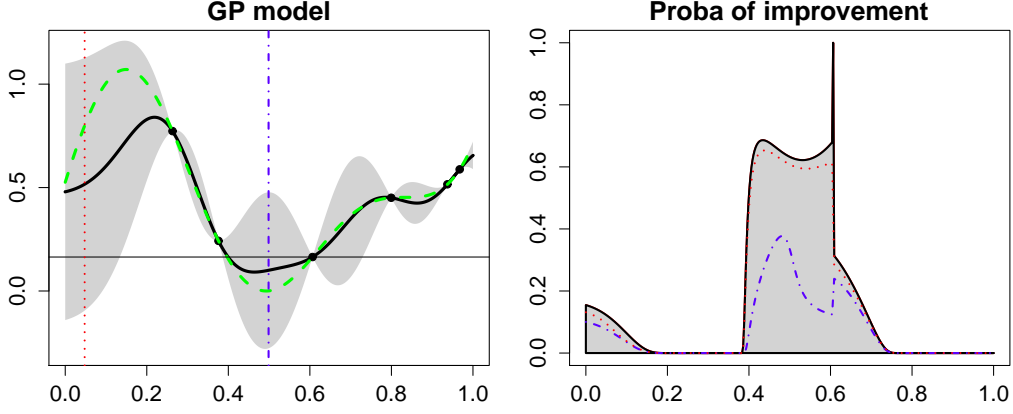}
\caption{Illustration of the effect of a new observation on the EV criterion. Left: actual objective function (dotted line), GP model (depicted by its mean in black plain line and 95\% confidence interval in grey) based on six observations (black circles). The horizontal line shows the current minimum; the vertical bars are placed at two candidate locations. Right: probability of improvement given by the current model and expected updated probability for each candidate. Adding a point at $x^+=0.5$ (mixed line) is likely to reduce substantially the probability, while adding a point at $x^+=0.05$ (dotted line) has little expected effect.}
	\label{fig:exSUR1}
\end{figure*}

Figure \ref{fig:exSUR1} illustrates the concept of reduction of volume of excursion on a toy example. A GP model is built on a six-point training set, from which the probability of improvement $p_6(\x,y^{\min}_6)$ is computed for every point in $\Xset=[0,1]$. We see that it can be intepreted as an indicator of the uncertainty we have about the location of the actual minimizer $\x^*=0.47$, as the model can only predict that $\x^*$ is likely to be between $0.4$ and $0.6$. Then, we consider two candidate points ($\x^+=0.2$ and $\x^+=0.5$) and compute, for each, the expected new probability (integrand in Eq.~(\ref{eq:def3})). We see that the probability is likely to remain mostly unchanged by adding the measurement at $\x^+=0.2$ (which is, indeed, a region with high response value), while it would be considerably reduced by adding a measurement at $\x^+=0.5$. In terms of volume of excursion set, we have $EEV(0.2) \approx ev_6$ (no reduction), while $EEV(0.5) \approx ev_6 /3$ (large reduction): the $EEV$ 
criterion clearly indicates $\x^+=0.5$ as a better sampling location.

\section{Multi-objective optimization by stepwise uncertainty reduction}\label{sec:multi}
\subsection{Volume of excursion behind the Pareto front}\label{sec:pareto}
Let $\y(\x) = \left(y^{(1)}(\x), \ldots, y^{(q)}(\x)\right)$ 
be the vector of objective functions to minimize. A point $\x$ dominates another point $\x'$ if 
$y^{(k)}(\x) \leq y^{(k)}(\x')$ for all $k$ in $\{1, \ldots, q\}$, which we denote by 
$\x' \prec \x$ in the following.
At step $n$, $\X_n = \{\x_1, \ldots, \x_n \}$ is the current experimental set and $\Y_n = \{\y(\x_1), \ldots, \y(\x_n)\}$ the corresponding set of measures. The non-dominated subset $\mathbf{X}_n^*$ of $\mathbf{X}_n$ constitutes the \textit{current Pareto front} (of size $m\leq n$).
In the objective space, the corresponding subset $\Y_n^*$ separates the regions dominated and not dominated by the experimental set. 

\begin{figure}[!ht]
	\centering
	\includegraphics[height=55mm]{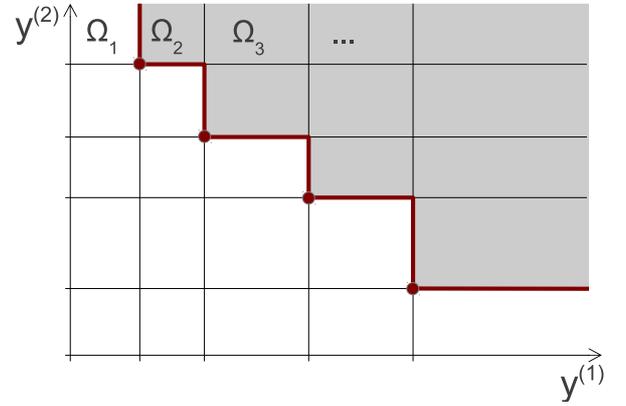}
	\caption{Example of Pareto front generated by four points (circles), and associated tesselation. The grey area corresponds to the dominated cells.}
	\label{fig:ParetoOnly}
\end{figure}

Then, we decompose the objective space plane using a tesselation $\{ \Omega_{i} \}_{i \in \{1, \ldots, I\}}$ of size $I=(m+1)^q$ ($\cup_{i \in I} \Omega_{i} = \mathbb{R}^q$ and $\cap_{i \in I} \Omega_{i} = \emptyset$), 
each cell being a hyperrectangle defined as:
\begin{equation}
 \Omega_{i} = \{ \y \in \mathbb{R}^q | \ykim \leq y^{(k)} < \ykip, k \in \{1, \ldots, q\} \}. \nonumber
\end{equation}
Each couple $(\ykim, \ykip)$ consists of two consecutive values of the vector $\left[-\infty, y^{(k)}(\x_1^*), \ldots, y^{(k)}(\x_m^*), +\infty \right]$.
An illustration is given in Figure \ref{fig:ParetoOnly}. 

A cell $\Omega_i$ \textit{dominates} another cell $\Omega_j$ ($\Omega_j \prec \Omega_i$) if any point in $\Omega_i$ dominates any point in $\Omega_j$, and it \textit{partially dominates} $\Omega_j$ if there exists a point in $\Omega_j$ that is dominated by any point in $\Omega_i$. Otherwise, we say that $\Omega_j$ is not dominated by $\Omega_i$ ($\Omega_j \not\prec \Omega_i$).

We denote by $I^*$ the indices of all the non-dominated cells at step $n$, that is, the cells that are not dominated by any point of $\mathbf{X}_n^*$. 
In two dimensions, the non-dominated cells are located in the bottom left half of the plane (Figure \ref{fig:ParetoOnly}).

\begin{figure}[!ht]
	\includegraphics[height=55mm]{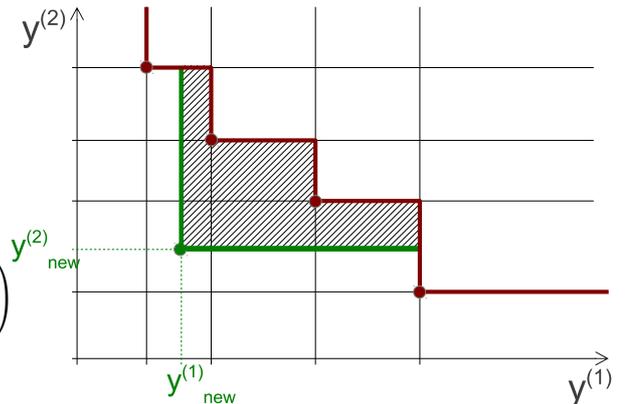}
	\caption{Example of Pareto front modification due to a new measurement. Two points are removed from the Pareto front while the new point is added. The hatched area represents the additional dominated region.}
	\label{fig:ParetoUpdate}
\end{figure}

Now, let us assume that GP models are fitted to each objective $y^{(k})$. At step $n$, the probability that $\Y(\x)$ belongs to the cell $\Omega_{i}$ is:
\begin{eqnarray}
 p_n^{i}(\x) &=& 
 \prob \left[ \Y(\x) \in \Omega_{i} \big| \An \right]  \nonumber  \\
 &=& \prod_{k=1}^q
 \Phi \left( \frac{ \ykip - m_n^{(k)}(\x)}{s_n^{(k)}(\x)}\right) - \Phi \left( \frac{\ykim - m_n^{(k)}(\x)}{s_n^{(k)}(\x)}\right) \nonumber \\
 &:=& \prod_{k=1}^q p_n^{i(k)}\nonumber
\end{eqnarray}
by pairwise independence of $Y^{(1)}, \ldots, Y^{(q)}$. 
The probability that $\x$ is not dominated by any point of $\X_n$ is then the probability that $\Y(\x)$ belongs to one of the non-dominated parts of the objective space. As the $\Omega_{i}$'s are disjoint, it is equal to:
\begin{equation}\label{eq:defND}
 \prob(\x \not\prec \X_n \big| \An) = \sum_{ i \in I^*} p_n^i(\x).
\end{equation}

Finally, the volume of the excursion sets behind the Pareto front is equal to the integral of this probability over $\Xset$:
\begin{eqnarray} \label{eq:paretovol}
 ev_n &=& \int_\Xset \prob(\x \not\prec \X_n \big| \An) d\x \nonumber \\
 &=&  \int_\Xset \sum_{ i \in I^*} p_n^i(\x)d\x = \sum_{ i \in I^*} \int_\Xset p_n^i(\x)d\x.
\end{eqnarray}
When $ev_n$ is high, a large proportion of the design space is likely to be better than the current Pareto set; inversely, when $\mathbf{X}_n^*$ approaches the actual Pareto set $\mathbb{X}^*$, the volume tends to zero. Hence, it defines naturally an uncertainty indicator for a SUR strategy.

\subsection{SUR criterion derivation}
Now, let us consider that a measurement $\y\+$ is performed at a point $\x\+$. Compared to step $n$, the volume $ev$ is modified by two means. First, the new measurement will modify the quantities $m_n^{(k)}(\x)$, $s_n^{(k)}(\x)$ ($k \in \{1, \ldots, q\}$), hence, the probabilities $p_n^i(\x)$. Second, if the new measurement is not dominated by the current Pareto set, it modifies the Pareto optimal front, as the new value $\y(\x\+)$ is added to $\Y^*$ and the values of $\Y^*$ dominated by $\y(\x\+)$ (if they exist) are removed. An example of such update is given in Figure \ref{fig:ParetoUpdate}.

Focusing on the probability that a point remains non-dominated (Eq.~(\ref{eq:defND})), accounting for the modifications of the models is relatively easy (that is, computing the quantity $p\+^i(.)$), but 
accounting for modifications in the Pareto front is complex, as both the number of elements and their values might change.
To address this issue, we consider that the updated probability 
$\prob(\x \not\prec \X\+ \big| \mathcal{A}\+)$ can be computed using the same sum as for $\prob(\x \not\prec \X_n \big| \An)$ (Eq.~(\ref{eq:defND})) by modifying its elements $p_n^{i}(\x)$:
\begin{equation}
 \prob(\x \not\prec \X\+ \big| \mathcal{A}\+) = \sum_{ i \in I^*} \tilde p\+^{i}(\x), \nonumber
\end{equation}
with 
 $$\tilde p\+^i(\x) = \prob \Big( \x \not\prec \X\+ \cap \Y(\x) \in \Omega_i \Big| \An, \y(\x\+)=\y\+ \Big)$$ and the $\Omega_i$'s defined using $\Y_n^*$ (not $\Y\+^*$). 

Seing from step $n$, the $\tilde P\+^{j}(\x)$ are random, as 
$$Y(\x\+)^{(k)} \sim \nlaw \left( m_n^{(k)}(\x\+), s_n^{(k)2}(\x\+) \right),$$ $\forall k \in \{1, \ldots, q\}$.

The expectation of the new volume is then:
\begin{eqnarray}
 EEV(\x\+) &=& \esp \left[ \int_{\Xset} \sum_{ j \in I^*} \tilde P\+^{j}(\x)d\x  \right] \nonumber \\
           &=& \sum_{ j \in I^*} \int_{\Xset} \esp \left[ \tilde P\+^{j}(\x) \right]d\x.\nonumber
\end{eqnarray}

This expression can be decomposed by conditioning on the range values of the new observation (using the fact that the $\Omega_i$'s are disjoint):
\begin{eqnarray}
 \esp \left[ \tilde P\+^{j}(\x) \right] &&\nonumber \\
 = \sum_{i \in I}
 \prob\+ \Big[&\x& \not\prec \X\+ \cap \Y(\x) \in \Omega_{j} \Big| \Y\+ \in \Omega_{i}  \Big] \nonumber \\
&\times& \prob\+ \Big[\Y\+ \in \Omega_{i} \Big] \nonumber \\
 := \sum_{i \in I} p_{ij}(\x), && \nonumber
\end{eqnarray}
where $\prob\+$ is the probability conditional on \\
$\Big( \An, \Y(\x\+)=\Y\+ \Big)$.

We first note from Proposition \ref{prop:4} that $$\prob\+ \left[\Y\+ \in \Omega_{i} \right]  = p_n^{i}(\x\+).$$
Then, leaving aside non-domination, the probability that $\Y(\x)$ belongs to $\Omega_j$ knowing that $\Y\+$ belongs to $\Omega_i$ is given by:
\begin{equation}
 \prob\+ \Big[\Y(\x) \in \Omega_{j} \Big| \Y\+ \in \Omega_{i}  \Big] \times p_n^{i}(\x\+) = \prod_{k=1}^q b_{ij}^{(k)}(\x), \nonumber
\end{equation}
with:
\begin{eqnarray}
 b_{ij}^{(k)}(\x) &=& \prob\+ \Big[\ykjm \leq Y^{(k)}(\x) < \ykjp \Big| \ykim \leq Y^{(k)}\+ < \ykip  \Big] \nonumber \\
 &\times& p_n^{i(k)}(\x\+), \nonumber
 \end{eqnarray}
\begin{eqnarray} 
 p_n^{i(k)}(\x\+) = \prob_n \left[ \ykim \leq Y^{(k)}\+ < \ykip  \right], \nonumber
\end{eqnarray}
by pairwise independence of $Y^{(1)}, \ldots, Y^{(q)}$. 
We show in Appendix \ref{sec:appendix3} that $b_{ij}^{(k)}(\x)$ can be expressed in closed form as:
\begin{eqnarray} 
 b_{ij}^{(k)}(\x) &= &
   \boldsymbol{\Phi}_{\rho}^{(k)} \left(\overline \ykip, \widetilde \ykjp \right) 
 - \boldsymbol{\Phi}_{\rho}^{(k)} \left(\overline \ykip, \widetilde \ykjm \right)\nonumber \\
 &-& \boldsymbol{\Phi}_{\rho}^{(k)} \left(\overline \ykim, \widetilde \ykjp \right)
 + \boldsymbol{\Phi}_{\rho}^{(k)} \left(\overline \ykim, \widetilde \ykjm \right)\nonumber
\end{eqnarray}
with the notations introduced in Section \ref{sec:update}.

Now, we define:
\begin{eqnarray} 
 d_{ij}^{(k)}(\x) &=& \prob\+ \Big[\ykjm \leq Y^{(k)}(\x) < \ykjp \cap Y\+^{(k)} \leq Y^{(k)}(\x) \nonumber \\
& \Big|& \ykim \leq Y^{(k)}\+ < \ykip  \Big] \times p_n^{i(k)}(\x\+),\nonumber
\end{eqnarray}
which is identical to $b_{ij}^{(k)}(\x)$ with the additional condition $Y^{(k)}\+ \leq Y^{(k)}(\x)$. 
This condition is met when the $k$-th component of the new observation dominates the $k$-th component of $\Y(\x)$. 
We have $\x \prec \x\+$ only if the condition $Y^{(k)}\+ \leq Y^{(k)}(\x)$ is met for all components, hence, with probability of occurence $\prod_{k=1}^p d_{ij}^{(k)}(\x)$. Three cases arise: 
\begin{itemize}
 \item $\ykim \geq \ykjp$: the component cannot be dominated, which implies $d_{ij}^{(k)}(\x)=0$;
 \item $\ykip \leq \ykjm$: the component is always dominated, which implies $d_{ij}^{(k)}(\x)=b_{ij}^{(k)}(\x)$;
 \item $\ykip = \ykjp$ (and $\ykim = \ykjm$): $Y^{(k)}(\x)$ and $Y^{(k)}\+$ share the same interval of variation, and:
\end{itemize}
\begin{eqnarray} 
 d_{ij}^{(k)}(\x) &=& \prob\+ \Big[F\+^{(k)} \leq Y^{(k)}(\x) < \ykip \Big| \ykim \leq Y^{(k)}\+ < \ykip  \Big] \nonumber \\
 &\times& p_n^{i(k)}(\x\+),\nonumber
\end{eqnarray} 
which is equal (as shown in Appendix \ref{sec:appendix3}) to:
\begin{eqnarray}
 d_{ij}^{(k)}(\x) &=&
   \boldsymbol{\Phi}_{\rho}^{(k)}\left(\overline \ykip, \widetilde \ykjp \right) 
 - \boldsymbol{\Phi}_{\nu}^{(k)}\left(\overline \ykip, \eta^{(k)} \right) \nonumber \\
 &+& \boldsymbol{\Phi}_{\nu}^{(k)}\left(\overline \ykim, \eta^{(k)} \right)
 - \boldsymbol{\Phi}_{\rho}^{(k)}\left(\overline \ykim, \widetilde \ykjp \right). \nonumber
\end{eqnarray}

The probability of $\Y(\x)$ being non-dominated while in $\Omega_j$ (and $\Y\+$ being in $\Omega_i$) is then:
\begin{equation}
 p_{ij}(\x) = \prod_{k=1}^q b_{ij}^{(k)}(\x) - \prod_{k=1}^q d_{ij}^{(k)}(\x).\nonumber 
\end{equation}

If $\Omega_{j} \prec \Omega_{i}$, 
the new observation dominates any point in $\Omega_{j}$, hence $d_{ij}^{(k)}(\x)=b_{ij}^{(k)}(\x)$ for all $k$, which gives $p_{ij}(\x) = 0$. 
Inversely, if $\Omega_{j} \not \prec \Omega_{i}$,
the new observation cannot dominate any point in the cell $\Omega_{j}$.
We have $d_{ij}^{(k)}(\x)=0$ for at least one value of $k$, and
$p_{ij}(\x)$ is the probability that $\Y(\x)$ belongs to $\Omega_{j}$: $p_{ij}(\x) = \prod_{k=1}^q b_{ij}^{(k)}(\x)$.

Finally, for a given point $\x \in \Xset$, we compute the probability that it is non-dominated at step $n+1$ using:
\begin{equation}
 \prob\+(\x \not\prec \X\+) = \sum_{i \in I} \sum_{j \in I^*} p_{ij}(\x),\nonumber
\end{equation}
and the SUR criterion is:
\begin{equation}
 EEV(\x\+) = \sum_{i \in I} \sum_{j \in I^*} \int _\Xset p_{ij}(\x) d\x,\label{eq:finalCrit}
\end{equation}
with:
\begin{equation}\label{eq:pijfinal}
 p_{ij}(\x) = 
 \left\{ 
 \begin{array}{ll}
    0 & \text{if } \Omega_{j} \prec \Omega_{i} \\
    \prod_{k=1}^q b_{ij}^{(k)}(\x) & \text{if } \Omega_{j} \not\prec \Omega_{i} \\
    \prod_{k=1}^q b_{ij}^{(k)}(\x) - \prod_{k=1}^q d_{ij}^{(k)}(\x) & \text{otherwise} 
 \end{array} \right. 
\end{equation}
The first sum in Eq.~(\ref{eq:finalCrit}) accounts for $\Y\+$ potentially being in any cell $\Omega_{i}$; the second sum accounts for $\Y(\x)$ potentially being in a non-dominated cell $\Omega_{j}$. 
\subsection{Computation}\label{sec:computation}
Evaluating the criterion as in Eq.~(\ref{eq:finalCrit}) is a non-trivial task; 
besides, a relatively fast computation is needed, as it may be embedded in an optimization loop to search for the best new observation (Eq.~(\ref{eq:SURstrategy})).
We provide here some technical solutions to ease its computation. Some of these issues have also been experienced with SUR criteria for inversion, as reported in \citet{chevalier2012fast,chevalier2013kriginv}.

Firstly, as no closed form exists for the integration over the design domain $\Xset$ in Eq.~(\ref{eq:finalCrit}), one may rely on Monte-Carlo integration, with approximations of the form:
\begin{equation}
 \int _\Xset p_{ij}(\x) d\x \approx \frac{1}{L} \sum_{l=1}^L w^{l} p_{ij}(\x^{l}),\nonumber
\end{equation}
where the $\x^{l}$'s and $w^{l}$'s are integration points and weights, respectively. 
One solution to alleviate the computational cost is to use a fixed set of integration points while searching for the best new observation. Then, many quantities that do not depend on $\x\+$ can be precalculated only once beforehand 
outside the optimization loop, as suggested in \citet{chevalier2012fast}.

Secondly, the criterion relies on the bivariate normal distribution, which also must be computed numerically. Very efficient programs can be found, such as the R package \texttt{pbivnorm} \citep{kenkel2011pbivnorm}, which makes this task relatively inexpensive. 

Thirdly, the tesselation used in the previous section has $I=(m+1)^q$ elements, with $I^*=I/2$ non-dominated elements, making the computation of the double sum in Eq.~(\ref{eq:finalCrit}) intensive. As detailed in Section \ref{sec:2D} for the two dimensional case, the number of elements can be very substantially reduced by grouping cells together. Note however that such decomposition may not be straightforward in high dimension.

Finally, as the optimization progresses, it is likely that the Pareto set grows, making the criterion more expensive to compute as more cells are to be considered. This problem is shared by all GP-based strategies, and some solutions have been proposed to filter the Pareto set and retain a small representative set \citep{wagner2010expected}. Such types of strategies may be applicable to our criterion, as some small cells would contribute to a very small part of the volume of excursion sets and could be neglected without introducing a critical error, and would reduce substantially the computational cost, especially when the number of observations is high. 
\subsection{Efficient formulas in the two-objective case}\label{sec:2D}
We consider here the two-objective case, for which the $EEV$ criterion can be expressed in a compact and computationally efficent way. 
With two objectives, the Pareto set can be ordered as follows (the first and second objective functions in ascending and descending order, respectively):
  $y_1^{(1)*} \leq \ldots \leq y_m^{(1)*}$ and 
  $y_1^{(2)*} \geq \ldots \geq y_m^{(2)*}$. 

\begin{figure*}[!ht]
	\centering
	\includegraphics[height=50mm]{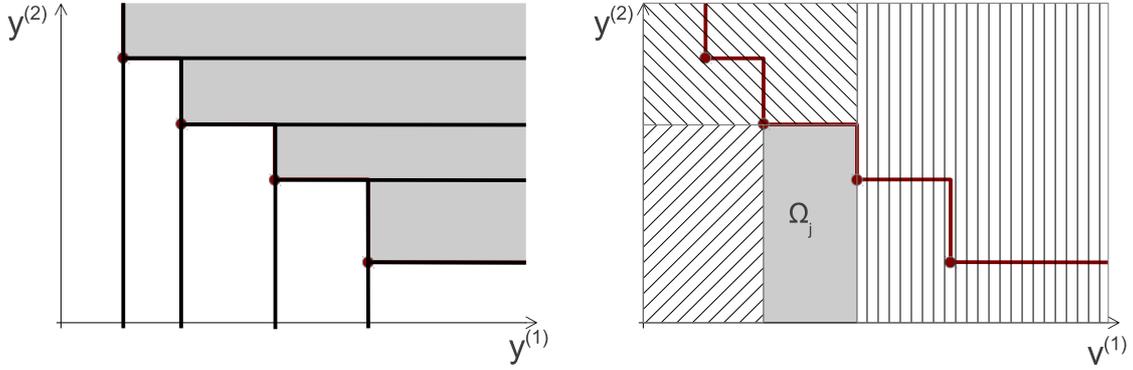}
	\caption{Left: the $m+1$ non-dominated cells (in white). Right: the four cases for $\Omega_i$ given $\Omega_j$: three are represented by the different hatched regions, the fourth corresponds to $\Omega_j=\Omega_i$.}
	\label{fig:Pareto2D}
\end{figure*}

The non-dominated part of the objective space can be divided in $m+1$ cells. Then, given a non-dominated cell $\Omega_j$, only four cases arise for $\Omega_i$ (the cell of the new observation), as shown in Figure \ref{fig:Pareto2D}, for which the quantities $p_{ij}(\x)$ need to be computed.

Hence, the criterion can be expressed as a sum of at most $(m+1)\times 4$ terms. As many terms can be factorized, we finally obtain:
\begin{equation}
 EEV(\x\+) = \sum_{j=0}^m \int_\Xset \alpha_j(\x),\nonumber
\end{equation}
with:
\begin{eqnarray}
 \alpha_0(\x) &=& \left[ \boldsymbol{\Phi}_{\rho}^{(1)} \left( \bar y_1^{(1)*},\tilde y_1^{(1)*} \right)  - \boldsymbol{\Phi}_{\nu}^{(1)} \left( \bar y_1^{(1)*},\eta^{(1)} \right) \right] \nonumber \\
 &\times& \left[ \Phi \left( \eta^{(2)}\right) - 1 \right]
 + \Phi \left( \tilde y_1^{(1)*}\right), \nonumber \\
 \alpha_j(\x) &=& 
 \big[ 
 \boldsymbol{\Phi}_{\rho}^{(1)} \left( \bar y_{j+1}^{(1)*}, \tilde y_{j+1}^{(1)*}  \right)
- \boldsymbol{\Phi}_{\nu}^{(1)} \left(\bar y_{j+1}^{(1)*}, \eta^{(1)}  \right) \nonumber \\
 &+&
 \boldsymbol{\Phi}_{\nu}^{(1)} \left(\bar y_{j}^{(1)*} , \eta^{(1)} \right) 
- \boldsymbol{\Phi}_{\rho}^{(1)} \left( \bar y_{j}^{(1)*}, \tilde y_{j}^{(1)*} \right)
 \big] \nonumber \\
 &\times&
 \left[  
  \boldsymbol{\Phi}_{\nu}^{(2)} \left(\bar y_{j}^{(2)*}, \eta^{(2)}  \right)
 -\boldsymbol{\Phi}_{\rho}^{(2)} \left(\bar y_{j}^{(2)*},  y_{j}^{(2)*}  \right)
 \right] \nonumber \\
 &+&
 \left[
 \Phi \left( \tilde y_{j+1}^{(1)*}\right) - 
 \Phi \left( \tilde y_j^{(1)*}\right)
 \right]
 \Phi \left( \tilde y_j^{(2)*}\right),
  \nonumber \\
  \forall j &\in& \{1, \ldots, m-1\}, \nonumber 
  \end{eqnarray}
  and:
  \begin{eqnarray}
 \alpha_m(\x) &=& 
 \big[  1 - \Phi \left( \eta^{(1)} \right) 
 + \boldsymbol{\Phi}_{\nu}^{(1)} \left( \bar y_m^{(1)*},\eta^{(1)} \right) 
 \nonumber \\
 &-& \boldsymbol{\Phi}_{\rho}^{(1)} \left( \bar y_m^{(1)*},\tilde y_m^{(1)*} \right) \big] \nonumber \\
 &\times&
 \left[  \boldsymbol{\Phi}_{\nu}^{(2)}   \left( \bar y_m^{(2)*},\eta^{(2)} \right)
 - \boldsymbol{\Phi}_{\rho}^{(2)} \left( \bar y_m^{(2)*},\tilde y_m^{(2)*} \right) \right] 
\nonumber \\
 &+&
 \left[ 1 - \Phi \left(\bar y_m^{(1)*} \right) \right]  \Phi \left( \bar y_m^{(2)*}  \right).
  \nonumber
\end{eqnarray}
Calculations are not detailed, as they are straightforward developments of Eq.~(\ref{eq:pijfinal}). 
The two extremal terms ($j=0$ and $j=m+1$) correspond to special cases of $\Omega_j$ (first and last cells in Figure \ref{fig:Pareto2D}, right).
\section{Numerical experiments}\label{sec:illustration}
\subsection{One-dimensional, bi-objective problem}
In this section, we apply the method to the following bi-objective problem: $F^{(1)}$ and $F^{(2)}$ are independent realizations of one-dimensional GPs, indexed by a 300-point regular grid on $[0,1]$, with a stationary Matern covariance with regularity parameter $\nu=3/2$ \citep[chapter 4]{rasmussen2006gaussian}. The variance and range parameters are taken as one and $1/5$, respectively. 

\begin{figure*}[!ht]
	\centering
	\includegraphics[height=70mm]{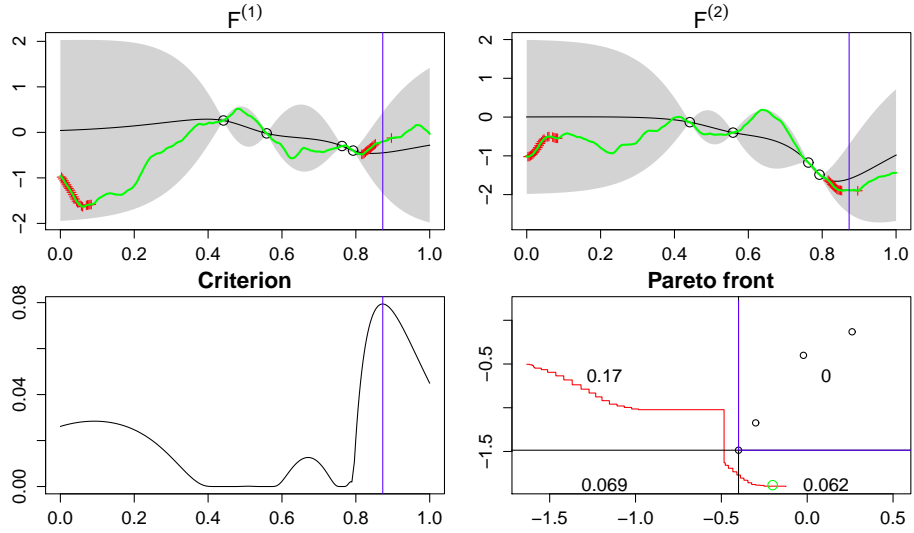}
\caption{Top graphs: initial models (same representation as Figure \ref{fig:exSUR1}); the actual Pareto-optimal points are represented by red crosses. Bottom right: observations (black circles) represented in the objective space, actual Pareto front (red) and current front (blue). Bottom left: criterion value as a function of $\x$. The vertical bars show the new observation location; the green circle is the new observation.}
	\label{fig:exmultiKrigInit}
\end{figure*}

Now, two GP models are built based on four randomly chosen observations. The covariance function is considered as known. Figure \ref{fig:exmultiKrigInit} shows the initial models and Pareto front. Here, a single point dominates the three others. After building the tesselation as described in Section \ref{sec:pareto}, we compute the volume of the excursion sets corresponding to each cell (Eq.~(\ref{eq:paretovol})). As there are only four observations, the probability to belong to a non-dominated cell is relatively high (Figure \ref{fig:exmultiKrigInit}, bottom right). Then, the criterion is computed for each point in the grid (Figure \ref{fig:exmultiKrigInit}, bottom right). Its maximum is obtained in a region with high uncertainty and low expected values for the two functions.

\begin{figure*}[!ht]
	\centering
\includegraphics[width=\textwidth]{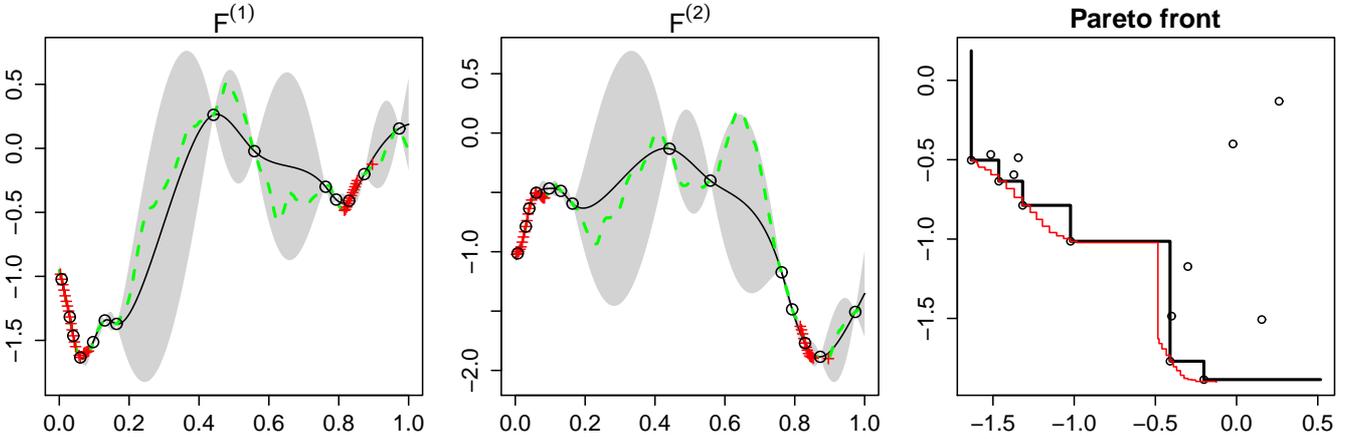}
\caption{Models and Pareto front after 10 iterations.}
	\label{fig:exmultiKrig10}
\end{figure*}

After 10 iterations (Figure \ref{fig:exmultiKrig10}), the Pareto front is well-approximated. The models are accurate in the optimal regions and have high prediction variances in the other regions, which indicates a good allocation of the computational resources.

\begin{figure*}[!ht]
	\centering
	\includegraphics[width=\textwidth]{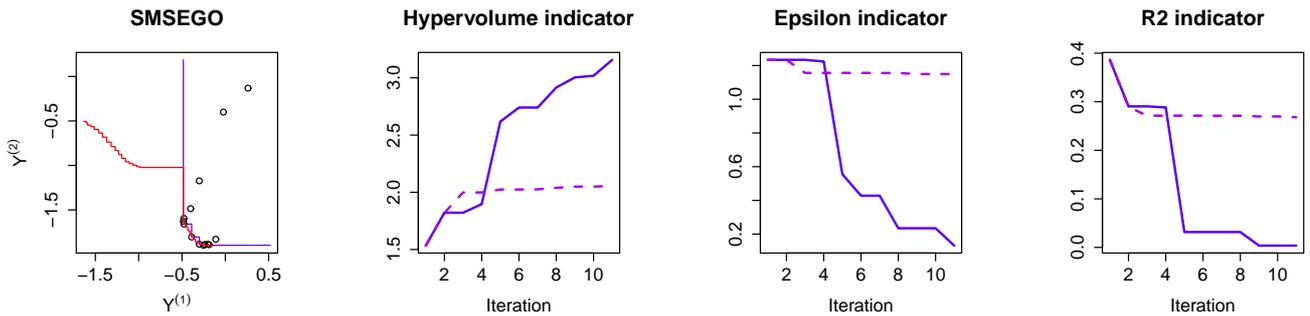}
        \caption{SMS-EGO Pareto front after 10 iterations (left) and performance comparison between SUR (plain line) and SMS-EGO (dotted line) on a one-dimensional problem. Hypervolume indicator: higher is better; other indicators should tend to zero.}
	\label{fig:SURvsSMSEGOdim2}
\end{figure*}

Next, we compare these results to a state-of-the-art method, called SMS-EGO \citep{ponweiser2008multiobjective}, which has been shown to outperform significantly non-GP based methods (such as NSGA-II), in particular when only a limited budget of evaluation is available. As measuring performances is non-trivial in multi-criteria optimization, we use a series of three indicators: hypervolume, epsilon and $R_2$ indicators \citep{zitzler2003performance,hansen1998evaluating}, all available in the R package EMOA \citep{mersmann2012emoa}. They provide different measures of distance to the actual Pareto set and coverage of the objective space. Results are reported in Figure \ref{fig:SURvsSMSEGOdim2}. The Pareto front obtained with SMS-EGO shows that the algorithm only detected one of the two Pareto optimal regions. As a consequence, the Pareto front is locally more accurate than the one obtained with the SUR strategy, but the indicators are much poorer.

\subsection{Six-dimensional, bi-objective problem}
\begin{figure*}[!ht]
	\centering
	\includegraphics[height=100mm]{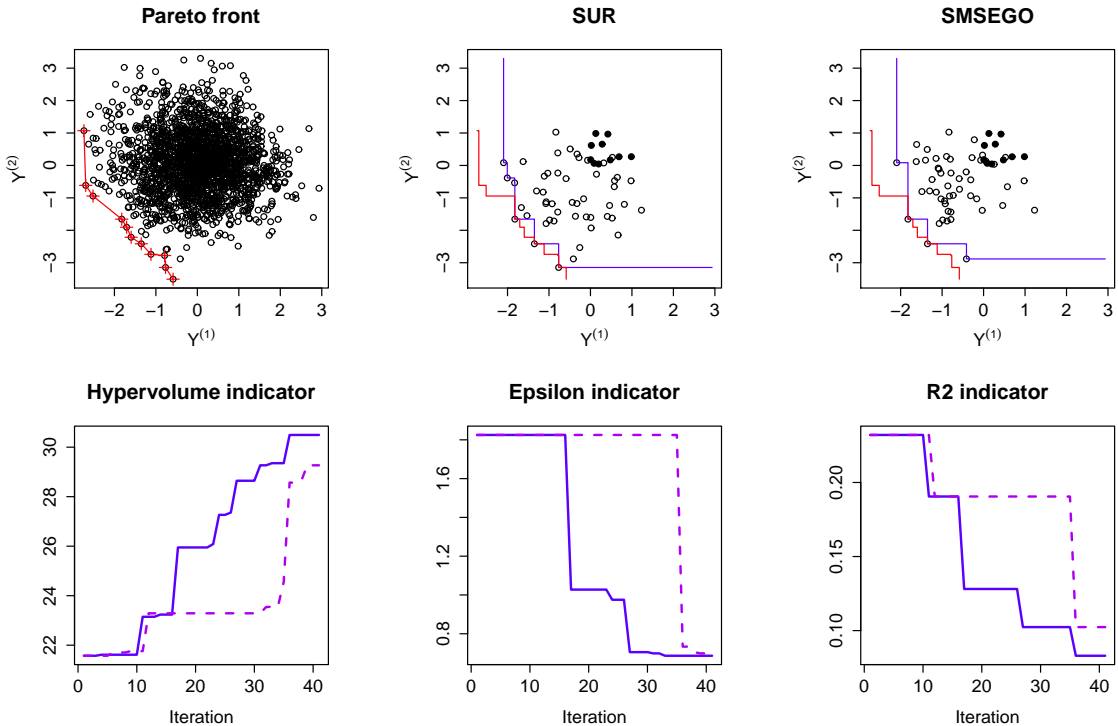}
        \caption{Performance comparison between SUR (plain line) and SMS-EGO (dotted line) on a six-dimensional problem. Top left: all 2000 points in the objective space; Pareto-optimal points are the red crosses. Top middle and right: Pareto fronts after 40 iterations.}
	\label{fig:SURvsSMSEGOdim6}
\end{figure*}

Here, the objectives functions are realizations of six-dimensional GPs indexed by a 2000-point Sobol sequence on $[0,1]^6$, with a stationary Matern covariance with regularity parameter $\nu=5/2$.
The variance and range parameters are taken as one and $ \sqrt{6}/6$, respectively. 
The initial experimental set consists of 10 points randomly chosen, and 40 points are added iteratively using the SUR and SMS-EGO strategies. The results are given in Figure \ref{fig:SURvsSMSEGOdim6}. Again, the SUR strategy shows better performances compared to SMS-EGO.

\section{Discussion}\label{sec:discussion}
We have proposed a new sequential sampling strategy, based on stepwise uncertainty reduction principles, for multi-objective optimization. Closed-form expressions were provided for the infill criterion. Numerical experiments showed promising performances of our strategy compared to a state-of-the-art method. We point here some strengths and weaknesses of our approach.

First of all, as it is based on Gaussian process modeling, it shares the limits inherents to the model. In particular, it is well-known that classical GP models cannot cope with large datasets ($>1000$) or high-dimensional spaces ($>100$). Most models also have restrictive conditions on the approximated function (typically, stationarity), and the strategy efficiency may be greatly penalized by important inadequations between the model hypothesis and the actual function characteristics. Using the proposed strategy on more complex GP models \citep{gramacy2008bayesian,banerjee2013efficient} may help mitigate these issues.

Secondly, we wish to emphasize here that the proposed method has a non-negligible computational cost, as (a) the criterion is evaluated by numerical integration and (b) it is embedded in an optimization loop. Hence, its use may be limited to simulators for which the time to compute an evaluation is much higher than the time to choose the next point to evaluate. However, one may note that the use of closed-form expressions, although relying on the bivariate normal CDF, avoid the need to use conditional simulations (as in \citet{villemonteix2009informational}) that would have made the method overly expensive.

On the other hand, moving away from the expected improvement paradigm allowed us to provide a method that does not necessitate any artificial ranking or trade-off between objective functions. It is also scale-invariant, which can be of great advantage when dealing with objectives of different nature. Finally, one advantage of the proposed strategy is that it considers progress in the design space rather than in the objective space, which corresponds to what practitioners are eventually interested in.

Possible extensions of this work are various. Accounting for the uncertainty due to the estimation of the model hyperparameters were left appart here; Bayesian approaches, in the fashion of \citet{kennedy2001bayesian} or \citet{gramacy2008bayesian} for instance, may help address this issue. Objective functions were considered as not correlated to ease calculations and allow the use of simple models. As objectives are likely to be negatively correlated in practice, accounting for it while keeping tractable criteria is an important question. Finally, the stepwise uncertainty reduction strategy may be easily adapted to other frameworks, such as constrained or noisy optimization.
\appendix
\section{Probabilities update}\label{sec:appendix2}
\subsection{Proof of Proposition \ref{prop:3}}\label{proof:3}
Using the model update equations (\ref{eq:updatedsn2}), we note first that:
\begin{equation}
  p\+(\mathbf{x},a) = \Phi \left[ \frac{a- m_n(\mathbf{x}) + \frac{c_n(\mathbf{x},\mathbf{x}\+)}{s_n^2(\mathbf{x}\+)} \left[ m_n\left(\mathbf{x}\+\right) - y\+  \right] }{s\+(\mathbf{x})} \right] \nonumber
\end{equation}
Now, let $\varphi\left( y\+ \right)$ be the PDF of $Y\+$ (conditional on $\An$). We have:
\begin{equation}
\begin{split}
 & q(\mathbf{x},  b, a) \nonumber \\
   =& \int_{-\infty}^b p\+(\mathbf{x}, a) d\varphi\left( y\+ \right) \nonumber \\
  =& \int_{-\infty}^b \Phi \left[ \frac{a- m_n(\mathbf{x}) + \frac{c_n(\mathbf{x},\mathbf{x}\+)}{s_n^2(\mathbf{x}\+)} \left[ m_n\left(\mathbf{x}\+\right) - y\+  \right]}{s\+(\mathbf{x})} \right] \nonumber \\
  & d\varphi\left( y\+ \right) \nonumber
  \end{split}
\end{equation}

As $Y\+ \sim \nlaw \left( m_n(\mathbf{x}\+), s_n^2(\mathbf{x}\+) \right)$, we can write (following \citet{chevalier2012fast}):
 $$Y\+ = m_n(\mathbf{x}\+) + s_n(\mathbf{x}\+)U$$ with $$U \sim \nlaw\left(0,1 \right),$$
which allows to simplify the previous equations to:
\begin{equation}\label{eq:prop3integform}
\begin{split}
 & q(\mathbf{x},  b, a)\\
  &= \int_{-\infty}^{\bar b}  \Phi \left[ 
 \frac{a - m_n(\mathbf{x})}{s\+(\mathbf{x})} - \left( \frac{c_n(\mathbf{x},\mathbf{x}\+)}{s_n(\mathbf{x}\+)s\+(\mathbf{x})} \right) u
 \right] d\varphi(u)  \\
 &= \int_{-\infty}^{\bar b} \Phi \left[ \hat a - \beta u \right] d\varphi(u), 
\end{split}
\end{equation}
with 
\begin{eqnarray}
 \beta &=& (c_n(\mathbf{x},\mathbf{x}\+))/(s_n(\mathbf{x}\+)s\+(\mathbf{x})), \nonumber \\
 \hat a &=& (a - m_n(\mathbf{x}))/s\+(\mathbf{x}) \text{ and}\nonumber \\
 \bar b &=& (b-m_n(\mathbf{x}\+))/s_n(\mathbf{x}\+). \nonumber
\end{eqnarray}
This quantity can be written as a bivariate Gaussian CDF. Indeed:
\begin{equation}
\begin{split}
&\int_{-\infty}^{\bar{b}}  \Phi \left[ \hat a - \beta u \right]  d\varphi(u)   \\
&  = \frac{1}{\sqrt{2 \pi}} \int_{-\infty}^{\bar{b}}  \Phi \left[ \hat a - \beta u \right] \exp \left( \frac{-u^2}{2} \right) du \nonumber\\
 &  = \frac{1}{2 \pi} \int_{-\infty}^{\bar{b}} \int_{-\infty}^{\hat a-\beta u} \exp \left[ - \frac{1}{2} \left(u^2+t^2 \right)\right] dt du \nonumber\\
 &  = \frac{1}{2 \pi} \int_{-\infty}^{\bar{b}} \int_{-\infty}^{\hat a} \exp \left[ - \frac{1}{2} \left(u^2+ \left[t-\beta u \right]^2 \right)\right] dt du \nonumber\\
&   = \frac{1}{2 \pi |\boldsymbol{\Sigma}_{\beta}|} \int_{-\infty}^{\bar{b}} \int_{-\infty}^{\hat a} \exp \left[ - \frac{1}{2} 
   \left[ \begin{matrix} u  &  t \end{matrix} \right]
   \boldsymbol{\Sigma}_{\beta} ^{-1}
   \left[ \begin{matrix} u  \\  t \end{matrix} \right]
   \right] dt du, \nonumber
\end{split}
\end{equation}
with $\boldsymbol{\Sigma}_{\beta} = 
\left[ \begin{matrix}
1  &  \beta\\
\beta  & 1 + \beta^2  \end{matrix} \right]$ 
(noting that $|\boldsymbol{\Sigma}_{\beta}|=1$), which is the standard form of the bivariate Gaussian CDF with zero mean and covariance matrix $\boldsymbol{\Sigma}_{\beta}$, hence:
\begin{equation}
 q(\mathbf{x},  b, a) = \boldsymbol{\Phi}_{\boldsymbol{\Sigma}_{\beta}}\left( \bar b, \hat a \right). \nonumber
\end{equation}
Finally, applying the normalization $\tilde a = \hat a / \sqrt{1 + \beta^2}$, we have: 
 $q(\mathbf{x},  b, a) = \boldsymbol{\Phi}_\rho\left( \bar b, \tilde a \right)$, 
with: $$\rho = \frac{\beta}{\sqrt{1 + \beta^2}} = \frac{c_n(\x,\x\+)}{s_n(\x\+)s_n(\x)}.$$
\subsection{Proof of Proposition \ref{prop:4}}\label{proof:4}
The result can be obtained directly from Proposition \ref{prop:3} with $b \rightarrow +\infty$. We have then: 
 $q(\mathbf{x},  b, a) \rightarrow \Phi\left( \tilde a \right) = p_n(\x, a)$.
\subsection{Proof of Corollary \ref{corr:1}}\label{proof:corr1}
From Eq.~(\ref{eq:prop3integform}), we have directly:
\begin{eqnarray}
   r(\mathbf{x},  b, a) &=& \int_{\bar b}^{+\infty} \Phi \left[ \hat a - \beta u \right] d\varphi(u) \nonumber \\
                               &=& \int_{-\infty}^{-\bar b} \Phi \left[ \hat a + \beta u \right] d\varphi(u)\nonumber \\
                               &=&  \boldsymbol{\Phi}_{\boldsymbol{\Sigma}_{-\beta}}\left( -\bar b, \hat a \right) = \boldsymbol{\Phi}_{-\rho}\left( -\bar b, \tilde a \right). \nonumber
\end{eqnarray}
\subsection{Proof of Proposition \ref{prop:5}}\label{proof:5}
The steps of the proof are similar to those of Proposition 3.
Using the update equations (\ref{eq:updatedsn2}), we have first:
\begin{equation}
\begin{split}
  &\mathbb{P} \left( Y(\mathbf{x}) \leq y\+ | \An, y\+=y(\mathbf{x}\+) \right) \\
  &= \Phi \left[ \frac{y\+ - m\+(\mathbf{x})}{s\+(\mathbf{x})} \right] \nonumber \\
  &= \Phi \left[ \frac{-m_n(\mathbf{x}) + \frac{c_n(\mathbf{x},\mathbf{x}\+)m_n\left(\mathbf{x}\+\right)}{s_n^2(\mathbf{x}\+)} +  \left[1 - \frac{c_n(\mathbf{x},\mathbf{x}\+)}{s_n^2(\mathbf{x}\+)} \right] y\+}{s\+(\mathbf{x})} \right] \nonumber\\
  &= \Phi \left[ \frac{m_n(\mathbf{x}\+) - m_n(\mathbf{x})}{s\+(\mathbf{x})} - \left( \frac{c_n(\mathbf{x},\mathbf{x}\+) - s_n^2(\mathbf{x}\+)}{s_n(\mathbf{x}\+)s\+(\mathbf{x})} \right) u \right]. \nonumber
  \end{split}
\end{equation}

Now:
\begin{equation}
\begin{split}
 & h(\mathbf{x},b) \\
 & = \int_{-\infty}^b \mathbb{P} \left( Y(\mathbf{x}) \leq Y\+ | \An, Y\+=y\+ \right) d\varphi\left( y\+ \right) \nonumber\\ 
 & = \int_{-\infty}^{\bar b} \Phi \Big[ \frac{m_n(\mathbf{x}\+) - m_n(\mathbf{x})}{s\+(\mathbf{x})} \\
 & - \left( \frac{c_n(\mathbf{x},\mathbf{x}\+) - s_n^2(\mathbf{x}\+)}{s_n(\mathbf{x}\+)s\+(\mathbf{x})} \right) u \Big] d\varphi(u) \nonumber \\
 & = \int_{-\infty}^{\bar b} \Phi \left[ \mu - \tau u \right] d\varphi(u) \nonumber \\
 & = \boldsymbol{\Phi}_{\boldsymbol{\Sigma}_{\tau}} \left(\bar b, \mu \right), \nonumber
\end{split}
\end{equation}
as we get a form similar to Equation \ref{eq:prop3integform}, with $\boldsymbol{\Phi}_{\boldsymbol{\Sigma}_{\tau}}$ the CDF of the centered bigaussian with covariance 
$\boldsymbol{\Sigma}_{\tau} = 
\left[ \begin{matrix}
1  &  \tau\\
\tau  & 1 + \tau^2  \end{matrix} \right]$, 
\begin{eqnarray}
 \mu &=& (m_n(\mathbf{x}\+) - m_n(\mathbf{x}))/s\+(\mathbf{x})\nonumber \\
 \tau &=& (c_n(\mathbf{x},\mathbf{x}\+) - s_n^2(\mathbf{x}\+))/(s_n(\mathbf{x}\+)s\+(\mathbf{x})). \nonumber
\end{eqnarray}
Normalizing $\eta = \mu/\sqrt{1 + \tau^2}$ delivers the final result.
\section{$b_{ij}^{(k)}(\x)$ and $d_{ij}^{(k)}(\x)$ computation}\label{sec:appendix3}
Let $X$ and $Y$ be two dependent random variables, and $a$, $b$, $c$ and $d$ four real numbers. By direct application of Bayes formula, we have:
\begin{equation}
\begin{split}
&\mathbb{P} \left( a \leq X < b | c \leq Y <d\right)  \mathbb{P} \left(c \leq Y <d\right) \\
&=
 \mathbb{P} \left(Y <d\right) \times \left[ \mathbb{P} \left( X < b | Y <d\right) - \mathbb{P} \left( X \leq a | Y <d \right) \right] \nonumber \\
 &- 
 \mathbb{P} \left(Y \leq c\right) \times \left[ \mathbb{P} \left( X < b | Y \leq c\right) - \mathbb{P} \left( X \leq a | Y \leq c \right) \right] \nonumber \\
\\
&\mathbb{P} \left( Y \leq X < b | a \leq Y < b\right)  \mathbb{P} \left(a \leq Y < b\right) \\
&=
 \mathbb{P} \left(Y < b\right) \times \left[ \mathbb{P} \left( X < b | Y < b\right) - \mathbb{P} \left( X \leq Y | Y < b \right) \right] \nonumber \\
 &-
 \mathbb{P} \left(Y \leq a\right) \times \left[ \mathbb{P} \left( X < b | Y \leq a\right) - \mathbb{P} \left( X \leq Y | Y \leq a \right) \right] \nonumber
 \end{split}
\end{equation}

Now, by definition, $b_{ij}^{(k)}$ is of the form of the fist equation:
\begin{eqnarray}
 b_{ij}^{(k)}(\x) &:=& \prob\+ \left( \ykjm \leq Y^{(k)}(\x) < \ykjp | \ykim \leq Y\+^{(k)} < \ykip \right)   \nonumber \\
 &\times&   \prob_n \left[ \ykim \leq Y\+^{(k)} < \ykip   \right], \nonumber
\end{eqnarray}
hence write as the sum of four terms:
\begin{eqnarray}
 b_{ij}^{(k)}(\x) 
 &=&  p_n^{(k)}(\x\+,\ykip) 
 \Big( \prob\+ \big[ Y^{(k)}(\x) \leq \ykjp \big| Y\+^{(k)}  \leq \ykip \big] \nonumber \\
 &-& \prob\+ \big[ Y^{(k)}(\x) \leq \ykjm \big| Y\+^{(k)} \leq \ykip \big] \Big) \nonumber \\
 &-& p_n^{(k)}(\x\+,\ykim) \Big( \prob\+ \big[ Y^{(k)}(\x) < \ykjp \big| Y\+^{(k)} \leq \ykim \big] \nonumber \\
 &-& \prob\+ \big[ Y^{(k)}(\x) \leq \ykjm \big| Y\+^{(k)}  \leq \ykim \big] \Big)   \nonumber \\
 &=& q^{(k)}\left(\x, \ykip, \ykjp \right) - q^{(k)}\left(\x, \ykip, \ykjm\right)  \nonumber \\
 &-& q^{(k)}\left(\x, \ykim, \ykjp \right) + q^{(k)}\left(\x, \ykim, \ykjm \right), \nonumber
\end{eqnarray}
with $q^{(k)}(\x, b, a)$ given by Eq.~(\ref{eq:prop3}), thus:
\begin{eqnarray}
 b_{ij}^{(k)}(\x) &=& 
   \boldsymbol{\Phi}_{\rho}^{(k)}\left(\overline \ykip, \widetilde \ykjp \right) 
 - \boldsymbol{\Phi}_{\rho}^{(k)}\left(\overline \ykip, \widetilde \ykjm \right)\nonumber\\
 &&- \boldsymbol{\Phi}_{\rho}^{(k)}\left(\overline \ykim, \widetilde \ykjp \right)
 + \boldsymbol{\Phi}_{\rho}^{(k)}\left(\overline \ykim, \widetilde \ykjm \right).\nonumber
\end{eqnarray}

Similary, $d_{ij}^{(k)}$ is of the form of the second equation:
Starting with the definition:
\begin{eqnarray}
 d_{ij}^{(k)}(\x) &:=&  \prob\+ \left( Y\+^{(k)} \leq Y^{(k)}(\x) < \ykjp \big| \ykim \leq Y\+^{(k)} < \ykip \right)  \nonumber \\ 
 &\times&   \prob_n \left[ \ykim \leq Y\+^{(k)} < \ykip   \right],\nonumber
\end{eqnarray}
hence writes:
\begin{eqnarray}
 d_{ij}^{(k)}(\x) 
 &=&  p_n^{(k)}(\x\+,\ykip) 
 \Big( \prob\+ \big[ Y^{(k)}(\x) \leq \ykjp \big| Y\+^{(k)}  \leq \ykip \big] \nonumber \\
 &-& \prob\+ \big[ Y^{(k)}(\x) \leq Y\+^{(k)} \big| Y\+^{(k)} \leq \ykip \big] \Big) \nonumber \\
 &-& p_n^{(k)}(\x\+,\ykim) \Big( \prob\+ \big[ Y^{(k)}(\x) < \ykjp \big| Y\+^{(k)} \leq \ykim \big] \nonumber \\
 &-& \prob\+ \big[ Y^{(k)}(\x) \leq Y\+^{(k)} \big| Y\+^{(k)}  \leq \ykim \big] \Big)   \nonumber \\
 &=& q^{(k)}\left(\x, \ykip, \ykjp \right) 
 - h^{(k)}\left(\x, \ykip \right) \nonumber \\
 &-& h^{(k)}\left(\x, \ykim \right)
 + q^{(k)}\left(\x, \ykim, \ykjm \right), \nonumber
\end{eqnarray}
with $q^{(k)}(\x, b, a)$ given by  Eq.~(\ref{eq:prop3}) and $h^{(k)}\left(\x, b \right)$ given by  Eq.~(\ref{eq:prop5}), thus:
\begin{eqnarray}
 d_{ij}^{(k)}(\x) &=& 
   \boldsymbol{\Phi}_{\rho}^{(k)}\left(\overline \ykip, \widetilde \ykjp \right) 
 - \boldsymbol{\Phi}_{\nu}^{(k)}\left(\overline \ykip, \eta^{(k)} \right) \nonumber \\
 &+& \boldsymbol{\Phi}_{\nu}^{(k)}\left(\overline \ykim, \eta^{(k)} \right)
 - \boldsymbol{\Phi}_{\rho}^{(k)}\left(\overline \ykim, \widetilde \ykjp \right). \nonumber
\end{eqnarray}



\end{document}